\documentclass{amsart}

\usepackage{latexsym}
\usepackage{amssymb,amsmath,amsfonts}

\newtheorem{theorem}{Theorem}[section]
\newtheorem{lemma}[theorem]{Lemma}
\newtheorem{cor}[theorem]{Corollary}

\theoremstyle{definition}

\newtheorem{example}[theorem]{Example}

\theoremstyle{remark}
\newtheorem{rmk}[theorem]{Remark}

\numberwithin{equation}{section}

\numberwithin{equation}{section}

\begin{document}

\title[Variational inequality problem]
 {Solution of variational inequality problems  on fixed point sets of nonexpansive mappings using iterative methods.}
 \address{$^1$Department of Mathematics, Nnamdi Azikiwe University, P. M. B. 5025, Awka, Anambra State, Nigeria.}

 \email{ euofoedu@yahoo.com}
\author[E. U. Ofoedu]{ Eric U. Ofoedu$^{*,\;1}$}

\keywords{ Modulus of smoothness, generalized duality maps, accretive operators,\\ strongly accretive mappings, Strictly pseudocontractive mappings, uniformly G\^{a}teaux differentiable norm.\\{\indent 2000 {\it Mathematics Subject Classification}. 47H06,
47H09, 47J05, 47J25.\\$^*$The author undertook this work when he visited
 the Abdus Salam International Centre
for Theoretical Physics (ICTP), Trieste, Italy as a visiting fellow.}}

\begin{abstract}
In this paper, we introduce  new implicit and explicit iterative schemes which converge strongly to a unique solution  of  variational inequality problems for strongly accretive operators over a common fixed point set of finite family of nonexpansive mappings in $q$-uniformly smooth real Banach spaces. As an application, we introduce an iteration process which converges strongly to a solution of the variational inequality which is a common fixed point of finite family of strictly pseudocontractive mappings. Our theorems extend, generalize, improve and unify the corresponding results of Xu \cite{27} and Yamada \cite{Yamada} and that of a host of other authors. Our corollaries and our method of proof are of independent interest.
\end{abstract}
\maketitle

\section{Introduction}
\noindent Let $(E, \;\| .\|)$ be a real normed space. Let $S:= \{x\in E: \|x\| = 1\}$. The space $E$ is said to have a
{\it G{\^ a}teaux differentiable} norm if the limit
\begin{eqnarray*}
\lim\limits_{t\to 0}\frac{\|x+ty\| - \|x\|}{t}
\end{eqnarray*}
exists for each $x, y\in S$. $E$ is said to have a {\it uniformly
G{\^a}teaux differentiable} norm if for each $y\in S$ the limit is
attained uniformly for $x\in S$. The space $E$ is said to be \emph{uniformly smooth} if and only of for all $\epsilon>0,$ there exists $\delta >0$ such that for all $x,y\in E$ with $\|x\|=1$ and $\|y\|\le \delta$, the inequality $$\frac{\|x+y\|+\|x-y\|}{2}-1<\epsilon \|y\|$$ holds. It is well known that every uniformly smooth real Banach space is \emph{a reflexive real Banach space} and has uniformly G\^ateaux differentiable norm (see e.g., \cite{ChidumeM}).\\

 \noindent Let $E$ be a real normed linear space with dimension, $dim(E)\ge 2$. The \emph{modulus of smoothness} of $E$ is the function $\rho_E:[0,\infty)\to [0,\infty)$ defined by $$\rho_E(t)=\sup\Bigl\{\frac{\|x+y\|+\|x-y\|}{2}-1:\|x\|=1,\|y\|=t\Bigr\}.$$ In terms of the modulus of smoothness (see e.g. \cite{ChidumeM}), the space $E$ is called \emph{uniformly smooth} if and only if $\displaystyle\lim_{t\to 0^+}\frac{\rho_E(t)}{t}=0.$ $E$ is called \emph{$q$-uniformly smooth} if and only if there exists a constant $c>0$ such that $\rho_E(t)\le ct^q, \;t>0$. It is easy to see that for $1<q<+\infty,$ every $q$-uniformly smooth Banach space is uniformly smooth; and thus has a uniformly G\^ateaux differentiable norm.\\
 
 \noindent Let $E$ be a real normed linear space  with dual $E^*.$
We denote by $J_q$ \emph{the generalized duality mapping} from $E$  to $2^{E^* }$ defined by
$$ J_qx:=\{f^* \in E^* :\langle x,f^*\rangle =\|x\|^q,\|f^*\|=\|x\|^{q-1} \}, $$  where $
\langle .,. \rangle $ denotes the generalized duality pairing between members of $E$ and 
members of $E^*$. For $q=2,$ the mapping $J=J_2$ from $E$ to $2^{E^*}$ is called \emph{the normalized duality mapping}. It
is well known that if $ E $ is uniformly smooth or $E^*$ is strictly convex, then duality mapping is
single-valued; and if $E$ has a uniformly G{\^a}teaux differentiable norm then the duality mapping is 
\emph{norm-to-weak$^*$ uniformly continuous} on bounded subsets of $E$. 
 If $E=H$ is a Hilbert space then the duality mapping becomes the identity map of $H$ (see e.g., \cite{ChidumeM, tbook}).
 In the sequel, we shall denote the single-valued generalized  duality mapping by
 $j_q$ and the single valued normalized duality mapping by $j$.
\\

\noindent A real normed space $E$ with strictly convex dual is said to have a \emph{weakly sequentially continuous generalized duality mapping} $j_q$ if and only if for each sequence $\{x_n\}_{n\ge 1}$ in $X$ such that $\{x_n\}_{n\ge 1}$ converges weakly to $x^*$ in $X$, we have that $\{j_q(x_n)\}_{n\ge 1}$ converges in the weak$^*$ topology to $j_q(x^*)$. A real Banach space $E$ is said to satisfy \emph{Opial's condition} if for any sequence $\{x_n\}_{n\ge 1}$ in $E$ such that $\{x_n\}_{n\ge 1}$ converges weakly to $x^*$ in $E$, we have that $$\displaystyle \limsup_{n\to\infty}\|x_n-x^*\|<\limsup_{n\to\infty}\|x_n-y\|$$ for all $y\in E,\;y\ne x^*.$ By Theorem 1 of Gossez and Lami Dozo \cite{Dozo}, it is known that if $E$ admits weakly sequentially continuous duality mapping, then $E$ satisfies Opial's condition (see \cite{Dozo} for more details). Aside from Hilbert spaces, it was noted in \cite{Browder1} that the most significant class of Banach spaces  having a weakly sequentially continuous generalized duality mapping are the sequence spaces $\ell_q$ for $1<q<+\infty$ (see \cite{Browder2} and \cite{Browder3}, where it is also shown that $L_q(\mathbb{R})$ has no weakly sequentially continuous duality mapping for $q\ne 2$).\\

\noindent A mapping $f:E\to E$ is said to be a \emph{strict contraction} or simply \emph{a contraction} if and only if there exists $\gamma_0\in [0,1)$ such that for all $x,y\in E,$ $$\|f(x)-f(y)\|\le \gamma_0\|x-y\|.$$ A mapping $T:E\to E$ is called \emph{nonexpansive} if and only if for all $x,y\in E,$ $$\|Tx-Ty\|\le \|x-y\|.$$ A point $x\in E$ is called a \emph{fixed point} of an operator $T$ if and only if $Tx=x.$ The set of fixed points of an operator $T$ is denoted by $F(T),$ that is, $F(T):=\{x\in E:Tx=x\}.$\\ 

\noindent A mapping $A:E\to E$ is said to be \emph{accretive} if and if for all $x,y\in E$ there exists $j_q(x-y)\in J_q(x-y)$ such that 
$$\Bigl<Ax-Ay,j_q(x-y)\Bigr>\ge 0.$$ The operator $A:E\to E$ is said to be \emph{$\psi$-uniformly accretive} if there exists a strictly increasing continuous function $\psi:[0,+\infty)\to [0,+\infty)$ with $\psi(0)=0$ such that for all $x,y\in E$, there exists $j_q(x-y)\in J_q(x-y)$ such that $$\Bigl<Ax-Ay,j_q(x-y)\Bigr>\ge \psi(\|x-y\|).$$  $A:E\to E$ is called \emph{strongly accretive} if and only if there exists a constant $\eta >0$ and for 
all $x,y\in E$ there exist $j_q(x-y)\in J_q(x-y)$ such that $$\Bigl<Ax-Ay,j_q(x-y)\Bigr>\ge \eta\|x-y\|^q.$$  It is easy to see that every strongly accretive mapping is $\psi$-uniformly accretive with $\psi(\tau)=\eta \tau^q\;\forall\;\tau\in [0,+\infty).$ When $E=H$ is a Hilbert space, accretive, $\psi$-uniformly accretive, strongly accretive mappings coincide with \emph{monotone, $\psi$-uniformly monotone, strongly monotone}  mappings, respectively.\\

\noindent A bounded linear operator $A:H\to H$ is called a \emph{strongly positive} operator if and only if there exist a constant $k>0$ such that for all $x\in H,$ $$\Bigl<Ax,x\Bigr>\ge k\|x\|^2.$$ Thus, every strongly positive bounded linear operator on $H$ is strongly monotone.\\

\noindent Iterative approximation of fixed points and zeros of nonlinear operators has been studied extensively by many authors to solve nonlinear operator equations as well as variational inequality problems (see e.g., \cite{ HZ3}, \cite{HZ6}-\cite{BA5}, \cite{HZ10}-\cite{HZ13}, \cite{7, 8}). In particular, iterative approximation of fixed points of nonexpansive mappings is an important subject in nonlinear operator theory and its applications in image recovery and signal processing is well known (see e.g., \cite{eu1, eu2, eu3}). Most published results on nonexpansive mappings centered on the iterative approximation of fixed points of the nonexpansive mappings or approximation of a common fixed points of a given family of nonexpansive mappings.\\

\noindent Iterative methods for nonexpansive mappings are now also applicable in solving convex minimization problems (see, for example, \cite{27} and references therein).  Let $H$ be a real Hilbert space with inner product $\bigl< ., .\bigr>$. Let $C$ be a closed convex nonempty subset of $H$, let $T:C\to C$ be a nonexpansive mapping such that $F(T)\ne \emptyset.$ Given $u\in C$ and a real sequence $\{\alpha_n\}_{n\ge 1}$ in the interval $(0,1)$, starting with an arbitrary initial $x_0\in C,$ let a sequence $\{x_n\}_{n\ge 1}$ be defined by
\begin{eqnarray}\label{Halpern}
x_{n+1}=\alpha_{n+1}u+(1-\alpha_{n+1})Tx_n, n\ge 0.
\end{eqnarray} Under appropriate conditions on the iterative parameter $\{\alpha_n\}_{n\ge 1},$ it has been shown by Halpern \cite{Halpern}, Lions \cite{Lions}, Wittmann \cite{Wittmann} and Bauschke \cite{Bauschke} that $\{x_n\}_{n\ge 0}$ converges strongly to $P_{F(T)}u,$ the projection of $u$ to the fixed point set, $F(T)$ of $T$. This means  that the limit of the sequence $\{x_n\}_{n\ge 0}$ solves the following minimization problem: 
$${\rm find}\; x^*\in F(T)\; {\rm such\; that\;}\displaystyle  \|x^*-u\|= \min_{x\in F(T)}\|x-u\|.$$

\noindent H. K. Xu \cite{27} studied the following quadratic minimization problem: find $x^*\in F(T)$ such that 
\begin{eqnarray}\label{X}
\displaystyle \frac{1}{2}\Bigl<Ax^*,x^*\Bigr>-\Bigl<x^*,u\Bigr>=\min_{x\in F(T)}\Bigl(\frac{1}{2}\Bigl<Ax,x\Bigr>-\Bigl<x,u\Bigr>\Bigr),
\end{eqnarray} where $u\in H$ is fixed and $A:H\to H$ a bounded linear strongly positive operator. Let $C_1, C_2,...,C_N$ be $N$ closed convex subsets of a real Hilbert space $H$ having a nonempty intersection $C$. Suppose also that each $C_i$ is a fixed point set of nonexpansive mappings $T_i:H\to H, i=1,2,..., N$, Xu \cite{27} proved strong convergence of the iterative algorithm 

\begin{eqnarray}\label{HKXU}
x_0\in H, \;x_{n+1}=(I-\alpha_{n+1}A)T_{n+1}x_n+\alpha_{n+1}u,\;n\ge 0
\end{eqnarray} (where $T_n=T_{n\;{\rm mod\;}N}$ and the mod function takes values in $\{1,2,...,N\}$) to a unique solution of the quadratic minimization problem \eqref{X}.\\

\noindent Marino and Xu \cite{Marino}, proved that the iteration scheme given by 
\begin{eqnarray}\label{4}
x_0\in H,\;x_{n+1}=\alpha_n\gamma f(x_n)+(I-\alpha_n A)Tx_n,\;n\ge 0
\end{eqnarray} converges strongly to a unique solution $x^\prime\in F(T)$ of the variational inequality 
\begin{eqnarray}
\Bigl<(\gamma f-A)x^\prime, y-x^\prime\Bigr>\le 0 \;\forall\;y\in F(T),
\end{eqnarray} which is the optimality condition for the minimization problem $$\min_{x\in F(T)}\Bigl(\frac{1}{2}\Bigl<Ax,x\Bigr>-h(c)\Bigr),$$ where $h$ is a potential function for $\gamma f$ (that is, $h^\prime(x)=\gamma f(x)$ for all $x\in H$); provided $f:H\to H$ is a contraction, $T:H\to H$ is nonexpansive and the iterative parameter $\{\alpha_n\}_{n\ge 0}$ satisfies appropriate conditions. \\

\noindent In \cite{Yamada}, Yamada introduced the following hybrid iterative method 
\begin{eqnarray}\label{5}
x_0\in H, \;x_{n+1}=Tx_n-\mu \lambda_n A(Tx_n), \;n\ge 0,
\end{eqnarray} where $T$ is nonexpansive, $A$ is $L$-Lipschitzian and strongly monotone operator with constant $\eta>0$ and $0<\mu <\frac{2\eta}{L^2}.$ He proved that if $\{\lambda_n\}_{n\ge 1}$ satisfy appropriate conditions, then \eqref{5} converges strongly to a unique solution $x^\prime\in F(T)$ of the varational inequality 
$$\Bigl<Ax^\prime, y-x^\prime\Bigl>\ge 0\;\forall\;y\in F(T).$$ 

\noindent Recently, M. Tian \cite{Tian} introduced the following iterative method:
\begin{eqnarray}
x_0\in H, \;x_{n+1}=\alpha_{n}\gamma f(x_n)+(I-\mu \alpha_n A)Tx_n, \;n\ge 0.
\end{eqnarray} Tian \cite{Tian} proved that if $f:H\to H$ is a contraction, $A:H\to H$ is an $\eta$-strongly monotone mapping, $T:H\to H$ a nonexpansive mapping and the parameter $\{\alpha_n\}_{n\ge 1}$ satisfies appropriate conditions, then the sequence $\{x_n\}_{n\ge 1}$ converges strongly to a unique solution $x^\prime\in F(T)$ of the variational inequality $$\Bigl<(\gamma f-\mu A)x^\prime, y-x^\prime\Bigr>\le 0\;\forall\;y\in F(T).$$ 

In this paper, motivated by the results of the authors mentioned above, it is our aim to introduce new implicit and explicit iterative schemes which converge strongly to a unique  solution of variational inequality problem 
over a common fixed point set of finite family of nonexpansive mappings in $q$-uniformly smooth real Banach spaces. As an application, we introduce an iteration process which converges storngly to a solution of the variational inequality which is a common fixed point of family of strictly pseudocontractive mappings. Our theorems extend, generalize, improve and unify the corresponding results of  Xu \cite{27} and Yamada \cite{Yamada} and that of a host of other authors. Our corollaries and our method of proof are of independent interest.

\section{Preliminary}

\noindent Let $\mu$ be a bounded linear functional defined on $\ell_\infty$ satisfying $\|\mu\|=1=\mu(1).$ It is known that $\mu$ is a mean on $\mathbb{N}$ if and only if $$\inf\{a_n:n\in\mathbb{N}\}\le \mu(a_n)\le\sup\{a_n:n\in\mathbb{N}\}$$ for every $a=(a_1,a_2,a_3,...)\in \ell_\infty.$
In the sequel, we shall use $\mu_n(a_n)$ instead of $\mu(a).$ A mean $\mu$ on $\mathbb{N}$ is called \emph{a Banach limit} if $\mu_n(a_n)=\mu_n(a_{n+1})$ for every $a=(a_1,a_2,a_3,...)\in \ell_\infty.$ It is well known that if $\mu$ is a Banach limit, then $$\displaystyle\liminf_{n\to\infty}a_n\le \mu_n a_n\le \limsup_{n\to\infty}a_n$$ for all $a=(a_1,a_2,a_3,...)\in \ell_\infty.$ Furthermore,
 if $a=(a_1, a_2, ...),\;b=(b_1,b_2,...)\in \ell_\infty$ and $\displaystyle\lim_{n\to\infty}a_n=a^*$ $\Bigl(\displaystyle\lim_{n\to\infty}(a_n-b_n)=0\Bigr)$, then
 $\mu_n(a_n)=a^*$ $\Bigl({\rm repectively,}\;\mu_n(a_n)=\mu_n(b_n)\Bigr)$.\\

\noindent In what follows, we shall make use of the following Lemmas:

\begin{lemma}\label{lem1}
Let $E$ be a real normed linear space, then for $1<q<+\infty,$ the following inequality holds: $$\|x+y\|^q\le \|x\|^q+q\big<y,j_q(x+y)\big>\;\forall\;x,y\in E,\;\forall\;j_q(x+y)\in J_q(x+y).$$ 
\end{lemma}

\begin{lemma}\label{lem2}
(See e.g., \cite{b, 26, 27}) Let $\{\lambda_n\}_{n\ge 1}$ be a sequence of nonnegative real numbers satisfying  the condition $$\lambda_{n+1}\le (1-\alpha_n)\lambda_n+\sigma_n,\;n\ge 0,$$ where $\{\alpha_n\}_{n\ge 0}$ and $\{\sigma_n\}_{n\ge 0}$ are sequences of real numbers such that $\displaystyle \{\alpha_n\}_{n\ge 1}\subset [0,1],\;\sum_{n=1}^\infty\alpha_n=+\infty.$ Suppose that $\sigma_n=o(\alpha_n),\;n\ge 0$ (i.e., $\displaystyle \lim_{n\to\infty}\frac{\sigma_n}{\alpha_n}=0$) or $\displaystyle \sum_{n=1}^\infty |\sigma_n|<+\infty$ or $\displaystyle\limsup_{n\to\infty}\frac{\sigma_n}{\alpha_n}\le 0$, then $\lambda_n\to 0$ as $n\to\infty.$
\end{lemma}

\begin{lemma}\label{lemmaT}(Compare with Lemma 3 pg. 257 of Bruck \cite{4}) Let $C$ be a nonempty closed and convex subset
of a real strictly convex Banach space $E$. Let $\{T_i\}_{i\ge 1}$ be a sequence of nonself nonexpansive mappings $T_i:C\to E$
such that $F:=\underset{i=1}{\overset{\infty}\bigcap}F(T_i)\ne \emptyset.$ Let $\{\sigma_i\}\subset (0,1)$ be such that $\displaystyle \sum_{i=1}^\infty\sigma_i=1.$ Then the mapping $T:=\displaystyle\sum_{i=1}^\infty\sigma_iT_i:C\to E$ is well defined, nonexpansive and $\displaystyle F(T)=\underset{i=1}{\overset{\infty}\bigcap}F(T_i).$
\end{lemma}

\begin{lemma}\label{lemma4} (See \cite{a28}, p.202 Lemma 3). Let $E$ be a
strictly convex Banach space and $C$ be a closed convex subset of
$E$. Let $T_1$, $T_2$,...,$T_r$ be nonexpansive mappings of $C$ into
itself such that the set of common fixed points of $T_1$,
$T_2$,...,$T_r$ is nonempty. Let $S_1$, $S_2$,...,$S_r$ be mappings
of $C$ into itself given by $S_i=(1-\gamma_i)I+\gamma_iT_i$ for any
$0<\gamma_i<1$, $i=1, 2, ...,r$, where $I$ denotes the identity
mapping on $C$. Then $S_1$, $S_2$,...,$S_r$ satisfies the following:
$\underset{i=1}{\overset{r}\bigcap}F(S_i)=\underset{i=1}{\overset{r}\bigcap}F(T_i)$ 
and 
$\underset{i=1}{\overset{r}\bigcap}F(S_i)=F(S_{r}S_{r-1}...S_{1})
=F(S_{1}S_{r-1}...S_{2})=...=F(S_{r-1}...S_{1}S_{r}).$
\end{lemma}

\begin{lemma}\label{lem4}(See e.g. \cite{Xu})
Let $E$ be a $q$-uniformly smooth real Banach space for some $q>1,$ then there exists some positive constant $d_q$ such that 
\begin{eqnarray}\label{qu}
\;\;\;\;\|x+y\|^q\le \|x\|^q+q\big<y,j_q(x)\big>+d_q\|y\|^q\;\forall\;x,y\in E,\;\forall\;j_q(x)\in J_q(x).
\end{eqnarray}
\end{lemma} 
\noindent If $E$ is $L_q$ (or $\ell_q$) space, the constant $d_q$ in \eqref{qu} has been calculated. This is shown in the following lemma.

\begin{lemma}\label{ell}(See e.g. \cite{HH3}) Let $E$ be $L_q$ (or $\ell_q$) space ($1<q<+\infty$) and $x,y\in E,$
\begin{enumerate}
\item if $1<q<2,$ then 
\begin{eqnarray*}
\|x+y\|^q\le \|x\|^q+q\big<y,j_q(x)\big>+d_q\|y\|^q\;\forall\;x,y\in E,\;\forall\;j_q(x)\in J_q(x),
\end{eqnarray*} where $d_q=\frac{1+b_q^{q-1}}{(1+b_q)^{q-1}},$ and $b_q$ is the unique solution of the equation \\
$(q-2)b^{q-1}+(q-1)b^{q-2}-1=0,$ $0<b<1.$
\item if $2\le q<+\infty,$ then 
\begin{eqnarray*}
\|x+y\|^2\le \|x\|^2+2\big<y,j(x)\big>+(q-1)\|y\|^2\;\forall\;x,y\in E,\;\forall\;j(x)\in J(x).
\end{eqnarray*}
\end{enumerate}
\end{lemma}

\begin{lemma}\label{deli}(Lemma 2.2 of \cite{Deli}, p. 1411)
Let $C$ be a closed convex nonempty subset of a reflexive Banach space which satisfies Opial's condition
and suppose $T:C\to E$ is nonexpansive, then the mapping $I-T$ is demiclosed at zero, that is is $\{x_n\}_{n\ge 1}$ is a sequence in 
$C$ such that $x_n\rightharpoonup x^*$ and $x_n-Tx_n\to 0$ as $n\to\infty,$ then $x^*=Tx^*$

\end{lemma}
\section{Main results.}
\begin{lemma}\label{eric}
Let $E$ be a real normed space. Let $A:E\to E$ be a strongly accretive mapping with a constant $\eta >0$. Let $\lambda>0$, then  the mapping $\lambda A:E\to E$ is strongly accretive.
\end{lemma}

\noindent{\bf proof.} Observe that 
\begin{eqnarray}\label{anaeme}
&&\Bigl<\lambda Ax-\lambda Ay,j_q(x-y)\Bigr>=\lambda\Bigl<Ax-Ay,j_q(x-y)\Bigr>\nonumber\\
&&\ge \lambda \eta\|x-y\|^q.
\end{eqnarray} This completes the proof. $\Box$\\

\begin{lemma}\label{uwam}
Let $E$ be a real normed space. Let $A:E\to E$ be a strongly accretive mapping with a constant $\eta >0$. Let $T:E\to E$ be a nonexpansive mapping such that $F(T)\ne \emptyset$, $\lambda >0$ and $u\in E$ be a fixed vector. Suppose that the solution $x^\prime\in F(T)$ of the variational inequality $\Bigl<u-\lambda Ax^\prime, j_q(p-x^\prime)\Bigr>\le 0\;\forall\;p\in F(T)$ exists, then $x^\prime$ is unique.
\end{lemma}

\noindent{\bf Proof.} Suppose for contradiction that the variational inequality has two solutions in $F(T),$ say $x^\prime \ne y^\prime,$ then we have, in particular, that $$\Bigl<u-\lambda Ax^\prime, j_q(y^\prime-x^\prime)\Bigr>\le 0$$ and $$-\Bigl<u-\lambda Ay^\prime, j_q(y^\prime-x^\prime)\Bigr>\le 0.$$ Thus, adding these two inequalities and using \eqref{anaeme}, we obtain that 
$$\lambda\eta\|y^\prime-x^\prime\|^q\le \Bigl<\lambda Ay^\prime-\lambda A x^\prime, j_q(y^\prime-x^\prime)\Bigr>\le 0,$$ a contradiction. Hence, the variational inequality has a unique solution, provided the solution exists. This completes the proof. $\Box$\\

\noindent Let $1<q<+\infty$ and let $E$ be a $q$-uniformly smooth real Banach space. Let $A:E\to E$ be an $L$-Lipschitzian strongly accretive mapping with a constant $\eta >0$. Let $T:E\to E$ a nonexpansive mapping. Then for the map $(I-t\lambda A)T:E\to E$ (where $I$ is the identity map of $K$, $\lambda>0$ and $t\in (0,1)$), we obtain using Lemma \ref{lem4} that
\begin{eqnarray*}
\|(I-t\lambda A)Tx-(I-t\lambda A)Ty\|^q &=&\|Tx-Ty-t\lambda(ATx-ATy)\|^q\\
&\le& \|Tx-Ty\|^q-qt\lambda\big<ATx-ATy,j_q(Tx-Ty)\big>\\&&+ d_q(t\lambda)^q\|ATx-ATy\|^q\\
&\le &\|Tx-Ty\|^q-qt\eta\lambda\|Tx-Ty\|^q\\&&+d_q(tL\lambda)^q\|Tx-Ty\|^q\\
&\le&\Bigl(1-t\lambda(q\eta-d_qL^q\lambda^{q-1})\Bigr)\|Tx-Ty\|^q
\end{eqnarray*}So that if $\lambda$ is such that $0 < \lambda<\Bigl(\frac{q\eta}{L^q d_q}\Bigr)^{\frac{1}{q-1}},$ we have that $$0<1-t\lambda(q\eta-d_qL^q\lambda^{q-1})<1$$ and since $T$ is a nonexpansive mapping of $E$ into $E$ and $1<q<+\infty$, we obtain

\begin{eqnarray}\label{nfor}
\|(I-t\lambda A)Tx-(I-t\lambda A)Ty\|&\le &\Bigl(1-t\lambda(q\eta-d_qL^q\lambda^{q-1})\Bigr)^{\frac{1}{q}}\|Tx-Ty\|\nonumber\\
&\le&\Bigl(1-t\lambda(q\eta-d_qL^q\lambda^{q-1})\Bigr)^{\frac{1}{q}}\|x-y\|.
\end{eqnarray}

\begin{lemma}\label{lom}
Let $E$ be a $q$-uniformly smooth real Banach space. Let $A:E\to E$ be an $L$-Lipschitzian strongly accretive mapping with a constant $\eta >0$. Let $T:E\to E$ be a nonexpansive mapping such that $F(T)\ne \emptyset$. Let $u\in E$ be fixed, $0 < \lambda<\Bigl(\frac{q\eta}{L^q d_q}\Bigr)^{\frac{1}{q-1}}$ and $t\in \Bigl(0,\min\Bigl\{1,\frac{1}{q\eta\lambda-d_q(L\lambda)^q}\Bigr\}\Bigr)$, then there exists a unique $z_t\in E$ such that 
\begin{eqnarray}\label{yu}
z_t=tu+(I-t\lambda A)Tz_t.
\end{eqnarray}
\end{lemma}

\noindent {\bf Proof.} For each $t\in \Bigl(0,\min\Bigl\{1,\frac{1}{q\eta\lambda-d_q(L\lambda)^q}\Bigr\}\Bigr),$ define $S_t:E\to E$ by $S_tx=tu+(I-t\lambda A)Tx$ for all $x\in E.$ Then, 
\begin{eqnarray*}
\|S_tx-S_ty\|&=& \|(I-t\lambda A)Tx-(I-t\lambda A)Ty\|\nonumber\\
&\le& \Bigl(1-t\lambda(q\eta-d_q(L)^q\lambda^{q-1})\Bigr)^{\frac{1}{q}}\|x-y\|\nonumber\\
\end{eqnarray*} Thus, for all $t\in \Bigl(0,\min\Bigl\{1,\frac{1}{q\eta\lambda-d_q(L\lambda)^q}\Bigr\}\Bigr)$, we have that $S_t$ is a strict contraction on $E$. Hence, there exists a unique $z_t\in E$ which satisfies \eqref{yu}. This completes the proof. $\Box$

\begin{lemma}\label{anya}Let $E$, $A$, $u\in E$ and $T$ be as in Lemma \ref{lom}. Let $\{z_t\}$ satisfy \eqref{yu}, then 
\begin{enumerate}
\item The mapping from $\Bigl(0,\min\Bigl\{1,\frac{1}{q\eta\lambda-d_q(L\lambda)^q}\Bigr\}\Bigr)$ to $E$ given by $t\mapsto z_t$ is continuous.
\item $\{z_t\}$ is bounded;
\item $\displaystyle \lim_{t\to 0}\|z_t-Tz_t\|=0;$
\end{enumerate}
\end{lemma}

\noindent {\bf Proof.} 
\begin{enumerate}
\item Let $t_0\in \Bigl(0,\min\Bigl\{1,\frac{1}{q\eta\lambda-d_q(L\lambda)^q}\Bigr\}\Bigr)$ be arbitrary. It is enough to show that the mapping from $\Bigl(0,\min\Bigl\{1,\frac{1}{q\eta\lambda-d_q(L\lambda)^q}\Bigr\}\Bigr)$ to $E$ given by $t\mapsto z_t$ is continuous at $t_0$.
Now, for some constant $M_1>0$ and using Lemma \ref{lem1}, we have that 
\begin{eqnarray}\label{view}
\|z_t-z_{t_0}\|^q&=&\|(t-t_0)u -(t-t_0)\lambda ATz_t+(I-t_0\lambda A)Tz_t-(I-t_0\lambda A)Tz_{t_0}\|^q\nonumber\\
&\le& \|(I-t_0\lambda A)Tz_t-(I-t_0\lambda A)Tz_{t_0}\|^q\nonumber\\&&+q\Bigl<(t-t_0)u-(t-t_0)\lambda ATz_t,j_q(z_t-z_{t_0})\Bigr>\nonumber\\
&\le& \Bigl(1-t_0(q\eta\lambda-d_q(L\lambda)^q)\Bigr)\|z_t-z_{t_0}\|^q+|t-t_0|M_1\|z_t-z_{t_0}\|^{q-1}.
\end{eqnarray}Thus, we obtain from \eqref{view} that $$\|z_t-z_{t_0}\|\le \frac{M_1}{t_0(q\eta\lambda-d_q(L\lambda)^q)}\Bigl|t-t_0\Bigr|$$ and the result follows.\\

\item Let $p\in F(T)$. If $z_t=p$ for all $t\in\Bigl(0,\min\Bigl\{1,\frac{1}{q\eta\lambda-d_q(L\lambda)^q}\Bigr\}\Bigr)$, then we are done. Otherwise, let $t^*\in \Bigl(0,\min\Bigl\{1,\frac{1}{q\eta\lambda-d_q(L\lambda)^q}\Bigr\}\Bigr)$ be such that $z_t\ne p$ $\forall\;t\in (0,t^*];$ then  using \eqref{nfor}, \eqref{yu} and Lemma \ref{lem1}, we have that $\forall\;t\in (0,t^*],$
\begin{eqnarray}\label{yu1}
\|z_t-p\|^q&= & \|(I-t\lambda A)Tz_t-(I-\lambda A)p+t(u-\lambda Ap)\|^q\nonumber\\
&\le & \Bigl(1-t(q\eta\lambda-d_q(L\lambda)^q)\Bigr)\|z_t-p\|^q\nonumber\\&&
+qt\Bigl<u-\lambda Ap,j_q(z_t-p)\Bigr>\nonumber\\
&\le&\Bigl(1-t(q\eta\lambda-d_q(L\lambda)^q)\Bigr)\|z_t-p\|^q\nonumber\\&&
+qt\|u-\lambda Ap\|.\|z_t-p\|^{q-1}
\end{eqnarray} Thus, we obtain from \eqref{yu1} that 
\begin{eqnarray}\label{yu23}
\|z_t-p\|\le \Bigl(1-t(q\eta\lambda-d_q(L\lambda)^q)\Bigr)\|z_t-p\|
+qt\|u-\lambda Ap\|.
\end{eqnarray}

Inequality \eqref{yu23} therefore gives $$\|z_t-p\|\le \frac{q\|u-\mu Ap\|}{q\eta\lambda-d_q(L\lambda)^q}.$$
Hence, the path $\{z_t\}$ is bounded; and so is $\{ATz_t\}$.\\

\item Therefore, using \eqref{yu}, we have that for some constant $M_0>0,$
$$\|z_t-Tz_t\|\le t\|u-\lambda ATz_t\|\le tM_0\to 0\;{\rm as\;} t\to 0.$$
This completes the proof. $\Box$
\end{enumerate}

\begin{theorem}\label{anya1}
Let $E$, $A$, $u\in E$ and $T$ be as in Lemma \ref{lom}. Let $\{z_t\}$ satisfy \eqref{yu}, then $\{z_t\}$ converges strongly to some $x^\prime\in F(T)$ which is a solution of the variational inequality 
\begin{eqnarray}\label{u}
&&\Bigl<u-\lambda Ax^\prime, j(p-x^\prime)\Bigr>\le 0\;\forall\;p\in F(T).
\end{eqnarray}
\end{theorem}

\noindent {\bf Proof.} By Lemma \ref{uwam}, if solution  of \eqref{u} exists in $F(T),$ then it is unique. Let $p\in F(T),$ then using \eqref{yu}, we have that
\begin{eqnarray}\label{mmo}
z_t-p=t\big(u-\lambda Ap\big)+(I-t\lambda A)Tz_t-(I-t\lambda A)p.
\end{eqnarray}
\noindent Thus, using \eqref{nfor}, \eqref{mmo} and Lemma \ref{lem1},
\begin{eqnarray}\label{uw}
\|z_t-p\|^q&=& \|(I-t\lambda A)Tz_t-(I-t\lambda A)p\|^q\nonumber\\
&&+qt\Bigl<u-\lambda Ap,j_q(z_t-p)\Bigr>\nonumber\\
&\le& \Big(1-t(q\eta\lambda-d_q(\lambda L)^q\Bigr)\|z_t-p\|^q\nonumber\\
&&+qt\Bigl<u-\lambda Ap,j_q(z_t-p)\Bigr>.
\end{eqnarray} So, \eqref{uw} implies that 
\begin{eqnarray}\label{yassir}
\|z_t-p\|^q\le \frac{1}{q\eta\lambda-d_q(\lambda L)^q}\Bigl<u-\lambda Ap,j_q(z_t-p)\Bigr>.
\end{eqnarray}

\noindent From Lemma \ref{anya} we know that $\{z_t\}$ is bounded. Let $\{t_n\}_{n\ge 1}$ in 
$ \Bigl(0,\min\Bigl\{1,\frac{1}{q\eta\lambda-d_q(L\lambda)^q}\Bigr\}\Bigr)$ be such that $\displaystyle \lim_{n\to\infty}t_n=0$ and set $z_n=z_{t_n}.$ Then, defining $\Psi:E\to \mathbb{R}\cup\{+\infty\}$ by $\Psi(x)=\mu_n\|z_n-x\|^q$ for all $x\in E$, where $\mu_n$ is a Banach limit on $\ell_{\infty}$, we have that $\Psi$ is continuous, convex and $\displaystyle\lim_{\|x\|\to+\infty}\Psi(x)=+\infty.$ Thus, setting $$\displaystyle K=\Bigl\{y\in E:\Psi(y)=\min_{x\in E}\Psi(x)\Bigr\},$$ then $K$ is bounded closed convex and nonempty subset of $E$ and since by Lemma \ref{anya}, $\displaystyle\lim_{n\to\infty}\|z_n-Tz_n\|=0, $ we obtain using Lemma \ref{lem1} that for some constant $Q^*>0$ and for all $y\in K$,
\begin{eqnarray}
\Psi(Ty)&=&\mu_n\|z_n-Ty\|^q = \mu_n\|Tz_n-Ty+z_n-Tz_n\|^q\nonumber\\
&\le &\mu_n\Bigl(\|Tz_n-Ty\|^q+q\Bigl<z_n-Tz_n,j_q(z_n-Ty)\Bigr>\Bigr)\nonumber\\
&\le &\mu_n\Bigl(\|Tz_n-Ty\|^q+\|z_n-Tz_n\|Q^*\Bigr)\nonumber\\
&\le&\mu_n\|z_n-y\|^q+\mu_n\|z_n-Tz_n\|Q^*= \mu_n\|z_n-y\|^q=\Psi(y).
\end{eqnarray} Thus, $T(K)\subset K$; that is, $K$ is invariant under $T$. Since $E$ is uniformly smooth real Banach space and thus every bounded closed convex nonempty subset of $E$ has the fixed point property for nonexpansive mappings, then there exists $x^\prime \in K$ such that $Tx^\prime =x^\prime.$ Since $x^\prime$ is also a minimizer of $\mu$ over $E$, it follows that for arbitrary $x\in E$ and for all $\xi\in (0,1),$ $\Psi(x^\prime)\le \Psi\big(x^\prime +\xi(x-x^\prime)\big).$ Lemma \ref{lem1} gives 
\begin{eqnarray}\label{eru}
&&\|z_n-x^\prime-\xi (x-x^\prime)\|^q\le \|z_n-x^\prime\|^q-q\xi\Bigl<x-x^\prime, j_q\big(z_n-x^\prime-\xi(x-x^\prime)\big)\Bigr>.
\end{eqnarray} Inequality \eqref{eru} implies that 
\begin{eqnarray}
\mu_n\Bigl<x-x^\prime, j_q\big(z_n-x^\prime-\xi(x-x^\prime)\big)\Bigr>\le 0.
\end{eqnarray} Furthermore, since $E$ is uniformly smooth, the generalized duality mapping $j_q$ is norm-to-norm uniformly continuous on bounded subsets of $E$; and hence
$$\displaystyle \lim_{\xi\to 0}\Bigl(\Bigl<x-x^\prime, j_q(z_n-x^\prime)\Bigr>-\Bigl<x-x^\prime, j_q\big(z_n-x^\prime-\xi(x-x^\prime)\big)\Bigr>\Bigr)=0.$$
So, for all $\epsilon>0,$ there exists $\delta_\epsilon >0$ such that for all $\xi \in (0,\delta_\epsilon)$ and for all $n\in \mathbb{N},$ we have that 
$$\Bigl<x-x^\prime, j_q(z_n-x^\prime)\Bigr>-\Bigl<x-x^\prime, j_q\big(z_n-x^\prime-\xi(x-x^\prime)\big)\Bigr><\epsilon.$$ 
This implies that for all $\xi \in (0,\delta_\epsilon)$, $$\mu_n\Bigl<x-x^\prime, j_q(z_n-x^\prime)\Bigr>\le \epsilon+\mu_n\Bigl<x-x^\prime, j_q\big(z_n-x^\prime-\xi(x-x^\prime)\big)\Bigr> \le \epsilon$$ and since $\epsilon>0$ is arbitrary, we obtain that $$\mu_n\Bigl<x-x^\prime, j_q(z_n-x^\prime)\Bigr>\le 0\;\forall\;x\in E.$$ Hence, we have in particular, that for the fixed 
$u\in E$ and $2x^\prime-\lambda Ax^\prime\in E$ 
\begin{eqnarray}\label{O}
\mu_n\Bigl<u-x^\prime, j_q(z_n-x^\prime)\Bigr>\le 0
\end{eqnarray} and 
\begin{eqnarray}\label{O1}&&
\mu_n\Bigl<(2x^\prime-\lambda Ax^\prime)-x^\prime, j_q(z_n-x^\prime)\Bigr>=\mu_n\Bigl<x^\prime-\lambda Ax^\prime,j_q(z_n-x^\prime)\Bigr>\le 0.
\end{eqnarray}
Moreover,  since $\{t_n\}_{n\ge 1}$ is in 
$ \Bigl(0,\min\Bigl\{1,\frac{1}{q\eta\lambda-d_q(L\lambda)^q}\Bigr\}\Bigr)$ and $x^\prime\in F(T),$ we obtain from \eqref{yassir} that
\begin{eqnarray}\label{laugh}
\|z_n-x^\prime\|^q&\le& \frac{1}{q\eta\lambda-d_q(\lambda L)^q}\Bigl<u-\lambda Ax^\prime,j_q(z_n-x^\prime)\Bigr>\nonumber\\
&=&\frac{1}{\omega} \Bigl(\Bigl<u-x^\prime, j_q(z_n-x^\prime)\Bigr>+\Bigl<x^\prime-\lambda Ax^\prime,j_q(z_n-x^\prime)\Bigr>\Bigr),
\end{eqnarray}where $\omega = q\eta\lambda-d_q(\lambda L)^q.$ So, from \eqref{laugh} (using \eqref{O}, \eqref{O1} and the linearity of Banach limit) we obtain that

\begin{eqnarray}
\mu_n\|z_n-x^\prime\|^q&\le& 
\frac{1}{\omega}\mu_n\Bigl<u-x^\prime, j_q(z_n-x^\prime)\Bigr>+\nonumber\\&&+\frac{1}{\omega}\mu_n\Bigl<x^\prime-\lambda Ax^\prime,j_q(z_n-x^\prime)\Bigr>\le 0.
\end{eqnarray} So, $\mu_n\|z_n-x^\prime\|^q=0$ and this implies that there exists a subsequence $\{z_{n_i}\}_{i\ge 1}$ of $\{z_n\}_{n\ge 1}$ such that $z_{n_i}\to x^\prime$ as $i\to\infty.$ \\

\noindent We now show that $x^\prime$ is a solution of the variational inequality \eqref{u}. From \eqref{yu}, we obtain that 
 $$z_{n_i}=t_{n_i}u+(I-t_{n_i}\lambda A)Tz_{n_i}.$$ This gives 

\begin{eqnarray}\label{dinna}
\lambda Az_{n_i}-u=-\frac{1}{t_{n_i}}(z_{n_i}-Tz_{n_i})+\lambda(Az_{n_i}-ATx_{n_i}).
\end{eqnarray} Thus, for all $p\in F(T),$ we obtain from \eqref{dinna} that

\begin{eqnarray}\label{uke1}
\Bigl<\lambda Az_{n_i}-u, j_q(z_{n_i}-p)\Bigr>&=&-\frac{1}{t_{n_i}}\Bigl<z_{n_i}-Tz_{n_i},j_q(z_{n_i}-p)\Bigr>\nonumber\\
&&+\lambda\Bigl<Az_{n_i}-ATx_{n_i},j_q(z_{n_i}-p)\Bigr>
\end{eqnarray}
Since $T$ is nonexpansive, we have that $(I-T)$ is accretive. Thus, $$\Bigl<z_{n_i}-Tz_{n_i},j_q(z_{n_i}-p)\Bigr>=\Bigl<(I-T)z_{n_i}-(I-T)p,j_q(z_{n_i}-p)\Bigr>\ge 0.$$ So, \eqref{uke1} gives

\begin{eqnarray}\label{uke}
\Bigl<\lambda Az_{n_i}-u, j_q(z_{n_i}-p)\Bigr>&\le&\lambda\Bigl<Az_{n_i}-ATx_{n_i},j_q(z_{n_i}-p)\Bigr>
\end{eqnarray}

\noindent Using \eqref{uke}, we obtain (for some constant $M_2>0$) that 
\begin{eqnarray}\label{uwad}
\Bigl<\lambda Ax^\prime-u, j_q(x^\prime-p)\Bigr>&=& \Bigl<\lambda Az_{n_i}-u, j_q(z_{n_i}-p)\Bigr>\nonumber\\
&&-\Bigl<\lambda Az_{n_i}-\lambda Ax^\prime, j_q(z_{n_i}-p)\Bigr>\nonumber\\
&&+\Bigl<u-\lambda Ax^\prime, j_q(z_{n_i}-p)-j_q(x^\prime -p)\Bigr>\nonumber\\
&\le& M_2\bigl(\|Az_{n_i}-ATz_{n_i}\|+\|\lambda Az_{n_i}-\lambda A x^\prime\|\nonumber\\
&&+\|j_q(z_{n_i}-p)-j_q(x^\prime -p)\|\Bigr).
\end{eqnarray} Hence, since $A,\; T$ are continuous, $z_{n_i}\to x^\prime\in F(T)$ as $i\to\infty$ and the generalized duality mapping is norm-to-norm uniformly continuous on bounded subsets of $E$ (and thus norm-to-norm continuous on bounded subsets of $E$), we obtain from \eqref{uwad} that (as $i\to\infty$)
$$\Bigl<\lambda Ax^\prime-u, j_q(x^\prime-p)\Bigr>\le 0\;\forall\;p\in F(T).$$ So, $x^\prime\in F(T)$ is a solution \eqref{u}.\\

\noindent Finally, we show that $\{z_n\}_{n\ge 1}$ converges to $x^\prime.$ Suppose that there is another subsequence $\{z_{n_l}\}_{l\ge 1}$ of $\{z_n\}_{n\ge 1}$ such that $z_{n_l}\to z^\prime\in E$ as $l\to\infty$. Then by serial number (2) of Lemma \ref{anya}, we have that $z^\prime\in F(T)$. Similar argument (from \eqref{dinna} to \eqref{uwad} with $z_{n_i}$ replaced by $z_{n_l}$) shows that $$\Bigl<\lambda Az^\prime-u, j_q(z^\prime-p)\Bigr>\le 0\;\forall\;p\in F(T).$$ Uniqueness of solution of \eqref{u} shows that $x^\prime=z^\prime.$ Hence, $z_n=z_{t_n}\to x^\prime$ and $n\to\infty.$ Consequently, we obtain that the path $\{z_t\}$ given by \eqref{yu} converges strongly (as $t\to 0$) to unique $x^\prime\in F(T)$ which solves the variational inequality \eqref{u}. This completes the proof. $\Box$\\
 
 \begin{rmk}
 One may worry about how $\displaystyle \lim_{t\to 0}z_t=x^\prime.$ Actually, this follows from the fact that for any sequence $\{t_m\}_{m\ge 1}$ in $ \Bigl(0,\min\Bigl\{1,\frac{1}{q\eta\lambda-d_q(L\lambda)^q}\Bigr\}\Bigr)$ such that $t_m\to 0$ as $m\to\infty,$ the argument of the proof of Theorem \ref{anya1} shows that $\{z_{t_m}\}_{m\ge 1}=\{z_m\}_{m\ge 1}$ converges to $x^*\in F(T)$ which solves the variational inequality \eqref{u} and since solution of \eqref{u} is unique, then $x^*=x^\prime.$ Thus, the path $\{z_t\}$ has a unique accumulation point as $t\to 0.$
 \end{rmk}
 
 \noindent The following corollaries follow from our discussion so far.

\begin{cor}\label{anya3}
Let $E$ be a real $L_p$ or $(\ell_p)$ space. Let $A$ and $T$ be as in Lemma \ref{lom} 
and $\{z_t\}$ satisfy \eqref{yu}, then we obtain the same conclusion as in Theorem \ref{anya1}.
\end{cor}

\begin{rmk}
Since every Hilbert space is a 2-uniformly smooth Banach space, it follows from Lemma \ref{ell} that if $E=H$ is a Hilbert space, then $d_q=d_2=1.$ Furthermore, we recall that in a Hilbert space $H$, the duality mapping coincide with the identity mapping on $H.$ Thus, we have the following corollary.
\end{rmk}
 
\begin{cor}\label{Hilbert}
Let $H$ be a real Hilbert space, $T:H\to H$ a nonexpansive mapping such that $F(T)\ne \emptyset$ and $A:H\to H$ be an $L$-Lipschitzian strongly accretive mapping with a constant $\eta >0$. Let $u\in H$ be fixed. Suppose that $0<\lambda<\frac{2\eta}{L^2}$ and $t\in \Bigl(0,\min\Bigl\{1,\frac{1}{2\eta\lambda-(\lambda L)^2}\Bigr\}\Bigr),$ then there exists a unique $z_t\in H$ satisfying \eqref{yu}. Moreover, $\{z_t\}$ converges strongly (as $t\to 0$) to a unique solution $x^\prime\in F(T)$ of the variational inequality 
$\Bigl<u-\lambda Ax^\prime, p-x^\prime\Bigr>\le 0\;\forall\;p\in F(T).$
\end{cor}

\begin{cor}\label{tre}
Let $E$ be a strictly convex $q$-uniformly smooth real Banach space, let $A:E\to E$ be an $L$-Lipschitzian strongly accretive mapping with a constant $\eta >0$. Let $T_i:E\to E, i=1,2,...$ be a countable family of nonexpansive mappings such that $F:=\underset{i=1}{\overset{\infty}\bigcap}F(T_i)\ne \emptyset.$ Let $\displaystyle T:=\sum_{i=1}^\infty\sigma_iT_i$, where $\{\sigma_i\}_{i\ge 1}\subset (0,1)$ is such that $\displaystyle \sum_{i=1}^\infty \sigma_i=1.$ Let $u\in E$ be fixed. Suppose that the conditions of Lemma \ref{lom} are satisfied, then  there  exist a unique $z_t\in E$ satisfying $z_t=tu+(I-t\lambda A)Tz_t.$
Moreover, $\{z_t\}$ converges strongly (as $t\to 0$) to a unique solution $x^\prime\in F(T)$ of the variational inequality 
$\Bigl<(u-\lambda Ax^\prime, j(p-x^\prime)\Bigr>\le 0\;\forall\;p\in F.$
\end{cor}

\noindent {\bf Proof.} By Lemma \ref{lemmaT}, $\displaystyle T:=\sum_{i=1}^\infty\sigma_iT_i$ is well defined, nonexpansive and $F(T)=F.$ The rest follows as in the proof of Theorem \ref{anya1}. $\Box$

\begin{cor}\label{tre1}
Let $E$ be a strictly convex $q$-uniformly smooth real Banach space, let $A:E\to E$ be an $L$-Lipschitzian strongly accretive mapping with a constant $\eta >0$. Let $T_i:E\to E, i=1,2,..., m$ be a finite family of nonexpansive mappings such that $F:=\underset{i=1}{\overset{m}\bigcap}F(T_i)\ne \emptyset.$ Let $\displaystyle T:=\sum_{i=1}^m\sigma_iT_i$, where $\{\sigma_i\}_{i=1}^m\subset (0,1)$ is such that $\displaystyle \sum_{i=1}^m \sigma_i=1.$ Let $u\in E$ be fixed. Suppose that the conditions of Lemma \ref{lom} are satisfied, then  there  exist a unique $z_t\in E$ satisfying $z_t=tu+(I-t\lambda A)Tz_t.$
Moreover, $\{z_t\}$ converges strongly to a unique solution $x^\prime\in F(T)$ of the variational inequality 
$\Bigl<u-\lambda Ax^\prime, j(p-x^\prime)\Bigr>\le 0\;\forall\;p\in F$.
\end{cor}

\section{Strong convergence of explicit iteration scheme for finite family of nonexpansive mappings.}
\noindent In the sequel, we shall assume that  $0 < \lambda<\Bigl(\frac{q\eta}{L^q d_q}\Bigr)^{\frac{1}{q-1}},$ $\alpha_n\in \Bigl(0,\min\Bigl\{1,\frac{1}{q\eta\lambda-d_q(L\lambda)^q}\Bigr\}\Bigr)$  $\forall\;n\in \mathbb{N}$, $\displaystyle \lim_{n\to\infty}\alpha_n=0,\;\sum_{n=1}^\infty \alpha_n=\infty$ and $\displaystyle \lim_{n\to\infty}\frac{\alpha_{n+r}-\alpha_{n}}{\alpha_{n+r}}=0\Longleftrightarrow \lim_{n\to\infty}\frac{\alpha_{n}}{\alpha_{n+r}}=1$.

\begin{theorem}\label{ICTP1}
Let $E$ be a $q$-uniformly smooth real Banach space which admits weakly sequentially continuous generalized duality mapping, let $A:E\to E$ be an $L$-Lipschitzian strongly accretive mapping with a constant $\eta >0$. Let $T_i:E\to E, i=1,2,..., r$ be a finite family of nonexpansive mappings such that $\Omega:=\underset{i=1}{\overset{r}\bigcap}F(T_i)\ne \emptyset.$   Let $u\in E$ be fixed and $\{x_n\}_{n\ge 1}$ be a sequence in $E$ generated iteratively by 
\begin{eqnarray}\label{Mama}
x_0\in E, x_{n+1}=\alpha_{n+1}u+(I-\alpha_{n+1}\lambda A)T_{n+1}x_n, \;n\ge 0,
\end{eqnarray}where $T_n=T_{n\;{\rm mod\;}r}$ and mod function takes values in $\{1,2,...,r\}.$ Suppose that $\Omega= F(T_{r}T_{r-1}...T_{1})=F(T_{1}T_{r-1}...T_{2})=...=F(T_{r-1}...T_{1}T_{r})$, then, $\{x_n\}_{n\ge 0}$ converges strongly to a solution of the variational inequality 
\begin{eqnarray}\label{Mama1}
&&\Bigl<u-\lambda Ax^\prime, j_q(p-x^\prime)\Bigr>\le 0\;\forall\;p\in \Omega.
\end{eqnarray}.
\end{theorem}

\noindent{\bf Proof.} Since the mapping $G=T_1T_2...T_r:E\to E$ is nonexpansive, then following the method of proof of Theorem \ref{anya1}, we have that for all $t\in \Bigl(0,\min\Bigl\{1,\frac{1}{q\eta\lambda-d_q(L\lambda)^q}\Bigr\}\Bigr)$, there exists unique $z_t\in E$ which satisfies $z_t=tu+(I-t\lambda A)Gz_t$, moreover, $\{z_t\}$ converges to a unique solution $x^\prime\in F(G)=\Omega$ of \eqref{Mama1}. We now show that the explicit scheme \eqref{Mama} converges strongly to $x^\prime.$  We start by showing first that $\{x_n\}_{n\ge 0}$ is bounded. Now, let $p\in \Omega$ and set $$r:= \max\Bigl\{\|x_0-p\|, \frac{q\|u-\lambda Ap\|}{q\eta\lambda-d_q(L\lambda)^q}\Bigl\}.$$ We show by induction that \begin{eqnarray}\label{recovery}
\|x_n-p\|\le r\;\forall\;n\ge 0.
\end{eqnarray} Observe that for $n=0$ \eqref{recovery} clearly holds. Assume for $n>0$ that \eqref{recovery} is true. We show that \eqref{recovery} is also true for $n+1$. Suppose for contradiction that this does not hold, then $$\|x_{n+1}-p\|>r\ge \|x_n-p\|.$$ Thus, using \eqref{nfor}, \eqref{Mama} and Lemma \ref{lem1}, we obtain that 
\begin{eqnarray}\label{Ukwe}
\|x_{n+1}-p\|^q&=&\|\alpha_{n+1}u+(I-\alpha_{n+1}\lambda A)T_{n+1}x_n-p\|^q\nonumber\\
&=&\|(I-\alpha_{n+1}\lambda A)T_{n+1}x_n\nonumber\\&&-(I-\alpha_{n+1}\lambda A)p+\alpha_{n+1}(u-\lambda Ap)\|^q\nonumber\\
&\le &\|(I-\alpha_{n+1}\lambda A)T_{n+1}x_n-(I-\alpha_{n+1}\lambda A)p\|^q\nonumber\\
&&+q\alpha_{n+1}\Bigl<u-\lambda Ap, j_q(x_{n+1}-p)\Bigr>\nonumber\\
&\le&\Bigl(1-\alpha_{n+1}(q\eta\lambda-d_q(L\lambda)^q)\Bigr)\|x_n-p\|^q\nonumber\\
&&+q\alpha_{n+1}\|u-\lambda Ap\|.\|x_{n+1}-p\|^{q-1}\nonumber\\
&<&\Bigl(1-\alpha_{n+1}(q\eta\lambda-d_q(L\lambda)^q)\Bigr)\|x_{n+1}-p\|^q\nonumber\\
&&+q\alpha_{n+1}\|u-\lambda Ap\|.\|x_{n+1}-p\|^{q-1}.
\end{eqnarray} Inequality \eqref{Ukwe} implies that $$\|x_{n+1}-p\|< \frac{q\|u-\lambda Ap\|}{q\eta\lambda-d_q(L\lambda)^q},$$ a contradiction.
Hence, the sequence $\{x_n\}_{n\ge 1}$ is bounded. Consequently, \\
$\{T_{n+1}x_n\}_{n\ge 1}$ and $\{AT_{n+1}x_n\}_{n\ge 1}$ are also both bounded.\\

\noindent Furthermore, using \eqref{Mama},
\begin{eqnarray}\label{bianu}
x_{n+r}-x_n&=&(\alpha_{n+r}-\alpha_n)u+(I-\alpha_{n+r}\lambda A)T_{n+r}x_{n+r-1}-(I-\alpha_{n+r}\lambda A)T_nx_{n-1}\nonumber\\&&-(\alpha_{n+r}-\alpha_n)\lambda AT_nx_{n-1}.
\end{eqnarray} Thus, using the fact that $T_{n+r}=T_n$, we obtain from \eqref{bianu} using Lemma \ref{lem1} that for some constant $M_5>0,$
\begin{eqnarray}\label{bianu1}
\|x_{n+r}-x_n\|^q&\le& \Bigl(1-\alpha_{n+1}(q\eta\lambda-d_q(L\lambda)^q)\Bigr)\|x_{n+r-1}-x_{n-1}\|^q\nonumber\\&&+|\alpha_{n+r}-\alpha_n|M_5\nonumber\\
&=&\Bigl(1-\alpha_{n+r}(q\eta\lambda-d_q(L\lambda)^q)\Bigr)\|x_{n+r-1}-x_{n-1}\|^q\nonumber\\&&+M_5\alpha_{n+r}\frac{|\alpha_{n+r}-\alpha_n|}{\alpha_{n+r}}
\end{eqnarray} So, using \eqref{bianu1}, we obtain from Lemma \ref{lem2} that 
\begin{eqnarray}\label{uch}
\displaystyle \lim_{n\to\infty}\|x_{n+r}-x_n\|=0
\end{eqnarray}
\noindent Furthermore, from the recursion formula \eqref{Mama} and for some constant $M_6>0$, we obtain
$$\|x_{n+1}-T_{n+1}x_n\|=\alpha_{n+1}\|u-\lambda AT_{n+1}x_n\|\le \alpha_{n+1}M_6.$$ 
Thus,
\begin{eqnarray*}
\displaystyle \lim_{n\to\infty}\|x_{n+1}-T_{n+1}x_n\|=0.
\end{eqnarray*} Similar argument shows that 
\begin{eqnarray}\label{bu}
x_{n+r}-T_{n+r}x_{n+r-1}\to 0\;{\rm as}\;n\to\infty.
\end{eqnarray}Thus, we obtain using \eqref{bu} and the fact that $T_n$ is nonexpansive that 
\begin{eqnarray*}
&&x_{n+r}-T_{n+r}x_{n+r-1}\to 0\;{\rm as}\;n\to\infty\nonumber\\
&&T_{n+r}x_{n+r-1}-T_{n+r}T_{n+r-1}x_{n+r-2}\to 0\;{\rm as}\;n\to\infty\nonumber\\
&&T_{n+r}T_{n+r-1}x_{n+r-2}-T_{n+r}T_{n+r-1}T_{n+r-2}x_{n+r-3}\to 0\;{\rm as}\;n\to\infty\nonumber\\
&&\vdots\nonumber\\
&&T_{n+r}...T_{n+2}x_{n+1}-T_{n+r}...T_{n+2}T_{n+1}x_n\to 0\;{\rm as}\;n\to\infty.\nonumber\\
\end{eqnarray*} Adding up these yields 
\begin{eqnarray}\label{ukk}
x_{n+r}-T_{n+r}...T_{n+2}T_{n+1}x_n\to 0\;{\rm as}\;n\to\infty.
\end{eqnarray}
\noindent But, 
\begin{eqnarray}\label{ukk5}
\|x_n-T_{n+r}...T_{n+2}T_{n+1}x_n\|&\le& \|x_n-x_{n+r}\|\nonumber\\&&+\|x_{n+r}-T_{n+r}...T_{n+2}T_{n+1}x_n\|
\end{eqnarray} Thus, using \eqref{uch}, \eqref{ukk} and \eqref{ukk5}, we get that 
\begin{eqnarray}\label{iwu}
\displaystyle \lim_{n\to \infty}\|x_n-T_{n+r}...T_{n+2}T_{n+1}x_n\|=0.
\end{eqnarray}

\noindent Next, we show that $$\displaystyle \limsup_{n\to\infty}\Bigl<u-\lambda Ax^\prime,j_q(x_n-x^\prime)\Bigr>\le 0.$$
Let $\{x_{n_k}\}_{k\ge 1}$ be a subsequence of $\{x_n\}_{n\ge 1}$ such that $$\displaystyle \limsup_{n\to\infty}\Bigl<u-\lambda Ax^\prime,j_q(x_n-x^\prime)\Bigr>=\lim_{k\to\infty}\Bigl<u-\lambda Ax^\prime,j_q(x_{n_k}-x^\prime)\Bigr>.$$ Since $\{x_n\}_{n\ge 1}$ is bounded and since $E$ is a reflexive real Banach space, there exists a subsequence $\{x_{n_{k_m}}\}_{m\ge 1}$ of $\{x_{n_k}\}_{k\ge 1}$ such that $\{x_{n_{k_m}}\}_{m\ge 1}$ converges weakly to some $p^*\in E.$ Without loss of generality, we may assume that $n_{k_m}$ is such that $T_{n_{k_m}}=T_i$ for some $i\in \{1,2,...,r\}$, for all $m\ge 1$. It therefore follows from \eqref{iwu} that $$\displaystyle\lim_{m\to\infty}\|x_{n_{k_m}}-T_{i+r}...T_{i+2}T_{i+1}x_{n_{k_m}}\|=0.$$ So, by Lemma \ref{deli}, $p^*\in F(T_{i+r}...T_{i+2}T_{i+1})=\Omega.$ Thus, since $E$ has weakly sequential continuous generalized duality mapping, we have that 
\begin{eqnarray}
\displaystyle \limsup_{n\to\infty}\Bigl<u-\lambda Ax^\prime,j_q(x_n-x^\prime)\Bigr>&=&\lim_{k\to\infty}\Bigl<u-\lambda Ax^\prime,j_q(x_{n_k}-x^\prime)\Bigr>\nonumber\\
&=&\lim_{m\to\infty}\Bigl<u-\lambda Ax^\prime,j_q(x_{n_{k_m}}-x^\prime)\Bigr>\nonumber\\
&=&\Bigl<u-\lambda Ax^\prime,j_q(p^*-x^\prime)\Bigr>\le 0\nonumber\\
\end{eqnarray} Thus, setting $$\theta_n=\max\Bigl\{0, \Bigl<u-\lambda Ax^\prime,j_q(x_n-x^\prime)\Bigr>\Bigr\},$$ then it is easy to see that 
$\displaystyle \lim_{n\to\infty}\theta_n=0.$ Furthermore,  we obtain from the recursion formula \eqref{Mama} using Lemma \ref{lem1} that 
\begin{eqnarray*}
\|x_{n+1}-x^\prime\|^q&\le & \Bigl(1-\alpha_{n+1}(q\eta\lambda-d_q(\lambda L)^q)\Bigr)\|x_n-x^\prime\|^q\nonumber\\ &&+q\alpha_{n+1}\Bigl<u-\lambda Ax^\prime,j_q(x_{n+1}-x^\prime)\Bigr>\nonumber\\
&\le &\Bigl(1-\alpha_{n+1}(q\eta\lambda-d_q(\lambda L)^q)\Bigr)\|x_n-x^\prime\|^q+\delta_n\; \forall\;n\ge 0,
\end{eqnarray*}where $\delta_n=\alpha_{n+1}q(q\eta\lambda-d_q(\lambda L)^q)\Bigl[ \frac{\theta_{n+1}}{q\eta\lambda-d_q(\lambda L)^q}\Bigr],$ which is clear $o(\alpha_{n+1}).$ Hence, by Lemma \ref{lem2}, $\{x_n\}_{n\ge 1}$ converges strongly to $x^\prime\in \Omega$ which is a unique solution of \eqref{Mama1}. This completes the Proof. $\Box$\\

\begin{rmk} We note that in the corresponding result of Xu \cite{27}, the assumption $$\Omega= F(T_{r}T_{r-1}...T_{1})=F(T_{1}T_{r-1}...T_{2})=...=F(T_{r-1}...T_{1}T_{r})$$ was made. 
We now consider a situation where this condition  is dispensed with.
\end{rmk}

\begin{example}\label{piu}
Let $E$ be a strictly convex $q$-uniformly smooth real Banach space which admits weakly sequentially continuous generalized duality mapping, let $A:E\to E$ be an $L$-Lipschitzian strongly accretive mapping with a constant $\eta >0$. Let $T_i:E\to E, i=1,2,..., r$ be a finite family of nonexpansive mappings such that $\Omega:=\underset{i=1}{\overset{r}\bigcap}F(T_i)\ne \emptyset$  and $S_i=(1-\omega_i)I+\omega_iT_i, i=1,2,...,r.$ Let $u\in E$ be fixed and let $\{x_n\}_{n\ge 1}$ be a sequence in $E$ generated iteratively by 
\begin{eqnarray}\label{Chimdindu1}
x_0\in E, x_{n+1}=\alpha_{n+1}u+(I-\alpha_{n+1}\lambda A)S_{n+1}x_n, \;n\ge 0,
\end{eqnarray}where $S_n=S_{n\;{\rm mod\;}r}$ and mod function takes values in $\{1,2,...,r\},$ then $\{x_n\}_{n\ge 0}$ converges strongly to a solution of the variational inequality \eqref{Mama1}
\end{example}

\noindent {\bf Proof.} By Lemma \ref{lemma4}, $\underset{i=1}{\overset{r}\bigcap}F(S_i)=\underset{i=1}{\overset{r}\bigcap}F(T_i)$ 
and 
$\underset{i=1}{\overset{r}\bigcap}F(S_i)=F(S_{r}S_{r-1}...S_{1})
=F(S_{1}S_{r-1}...S_{2})=...=F(S_{r-1}...S_{1}S_{r}).$ The rest follows as in the proof of Theorem \ref{ICTP1}. $\Box$

\begin{cor}\label{emphasis}
Let $E$ be a $q$-uniformly smooth real Banach space which admits weakly sequentially continuous generalized duality mapping, let $A:E\to E$ be an $L$-Lipschitzian strongly accretive mapping with a constant $\eta >0$. Let $T:E\to E,$ be a nonexpansive mapping such that $F(T)\ne \emptyset.$ Let $u\in E$ be fixed and let $\{x_n\}_{n\ge 1}$ be a sequence in $E$ generated iteratively by 
\begin{eqnarray}
x_0\in E, x_{n+1}=\alpha_{n+1}u+(I-\alpha_{n+1}\lambda A)Tx_n, \;n\ge 0,
\end{eqnarray} then $\{x_n\}_{n\ge 0}$ converges strongly to a unique $x^\prime \in F(T)$ which is a solution of the variational inequality \eqref{u}
\end{cor}

\noindent {\bf Proof.} Follows as in the proof of Theorem \ref{ICTP1} using Theorem \ref{anya1}. $\Box$

\begin{rmk}
We remark that corollaries synonymous to Corollaries \ref{anya3}, \ref{Hilbert}, \ref{tre} and \ref{tre1} are obtainable in this section. But we must note, however, that though Theorem \ref{ICTP1}, Corollary \ref{piu} and corollary \ref{emphasis} hold in the sequence space $\ell^q,$ they do not hold in $L^q(\mathbb{R})$ for $1<q<+\infty,\; q\ne 2$ since $L^q(\mathbb{R}),\;q\ne 2$ do not possess weakly sequentially continuous duality mapping. 
\end{rmk}

\begin{rmk}
The addition of bounded error terms to our recursion
formulas leads to no further generalization. It is easy to see that if we replace $u$ by $\lambda u$ in both implicit and explicit iteration scheme studied in this paper, we obtain that our iteration schemes converge strongly to a unique solution $x^\prime\in \Omega$ of the variational inequality $\Bigl<u-Ax^\prime, j_q(p-x^\prime)\Bigr>\le 0\;\forall\;p\in \Omega,$ where $\Omega$ is the set of common fixed points of finite family of nonexpansive mappings.
\end{rmk}

\section{Applications\\Convergence Theorem for families of strictly pseudocontractive mappings.}

\noindent Let $E$ be a normed space. A mapping $T:E\to E$ is called \emph{$k$-strictly pseudocontractive} if and only if there exists a real constant $k>0$ such that for all $x,y\in D(T)$ there exists $j(x-y)\in J(x-y)$ such that
\begin{eqnarray}\label{strictly}
\big<Tx-Ty,j(x-y)\big>\le \|x-y\|^2-k \|x-y-(Tx-Ty)\|^2.
\end{eqnarray}Without loss of generality we may assume that $k\in (0,1).$ If $I$ denotes the identity operator, then \eqref{strictly} can be re-written as
\begin{eqnarray}\label{sp}
\big<(I-T)x-(I-T)y,j(x-y)\big>\ge k \|(I-T)x-(I-T)y)\|^2.
\end{eqnarray}In Hilbert spaces, \eqref{strictly} (or equivalently \eqref{sp}) is equivalent to the inequality
$$\|Tx-Ty\|^2\le \|x-y\|^2+\beta\|(I-T)x-(I-T)y\|^2,\;{\rm where\;}\beta=(1-k)<1.$$

\noindent It was shown in \cite{Osilike} that if $T$ is $k$-strictly pseudocontractive, then the following inequality holds
\begin{eqnarray}\label{OsA}
\big<(I-T)x-(I-T)y,j_q(x-y)\big>\ge k^{q-1} \|(I-T)x-(I-T)y)\|^q.
\end{eqnarray}

\noindent Thus, if $E$ is a $q$-uniformly smooth real Banach space; and $T:E\to E$ is a $k$-strictly pseudocontractive mapping, then for the map $T_a:=(1-a)I+aT):E\to E$ (where $I$ is the identity map of $E$ and $a>0$), we obtain by Lemma \ref{lem4} using \eqref{OsA} that:
\begin{eqnarray*}
\|T_a x-T_a y\|^q &=&\|x-y-a\Bigl((I-T)x-(I-T)y\Bigr)\|^q\\
&\le& \|x-y\|^q
-q a\big<(I-T)x-(I-T)y,j_q(x-y)\big>\nonumber\\
&&+d_q a^q\|(I-T)x-(I-T)y\|^q\\
&\le &\|x-y\|^q-a(k^{q-1}q-d_q a^{q-1})\|Ax-Ay\|^q,
\end{eqnarray*} where $A=(I-T).$ If $a$ is such that $0 < a<\Bigl(\frac{qk^{q-1}}{d_q}\Bigr)^{\frac{1}{q-1}},$ we have that the mapping $T_a$ is nonexpansive. It is also easy to see that the fixed point set of $T_a$ and that of $T$ coincide.\\

\noindent Thus, we have the following theorem

\begin{theorem}
Let $E$ be a $q$-uniformly smooth real Banach space which admits weakly sequential continuous generalized duality mapping, let $A:E\to E$ be an $L$-Lipschitzian strongly accretive mapping with a constant $\eta >0$. Let $T_i:E\to E, i=1,2,..., r$ be a finite family of $k$-striclty pseudocontractive mappings such that $\Omega^*:=\underset{i=1}{\overset{r}\bigcap}F(T_i)\ne \emptyset.$ Let $\{a_i\}_{i=1}^r$ be such that
 $0 < a_i<\Bigl(\frac{q k^{q-1}}{d_q}\Bigr)^{\frac{1}{q-1}}, i=1,2,...,r$ and define $T_{a_i}=(1-a_i)I+a_i T_i.$  Let $u\in E$ be fixed and $\{x_n\}_{n\ge 1}$ be a sequence in $E$ generated iteratively by 
\begin{eqnarray}
x_0\in E, x_{n+1}=\alpha_{n+1}u+(I-\alpha_{n+1}\lambda A)T_{a_{n+1}}x_n, \;n\ge 0,
\end{eqnarray}where $T_{a_n}=T_{a_{n\;{\rm mod\;}r}}$ and mod function takes values in $\{1,2,...,r\}.$ Suppose that $\Omega^* = F(T_{a_r}T_{a_{r-1}}...T_{a_1})=F(T_{a_1}T_{a_{r-1}}...T_{a_2})=...=F(T_{a_{r-1}}...T_{a_1}T_{a_r})$, then, $\{x_n\}_{n\ge 0}$ converges strongly to a solution of the variational inequality 
\begin{eqnarray}\label{Ma}
&&\Bigl<u-\lambda Ax^\prime, j_q(p-x^\prime)\Bigr>\le 0\;\forall\;p\in \Omega^*.
\end{eqnarray}.
\end{theorem}

\begin{cor}
Let $E$ be a strictly convex $q$-uniformly smooth real Banach space which admits weakly sequentially continuous generalized duality mapping, let $A:E\to E$ be an $L$-Lipschitzian strongly accretive mapping with a constant $\eta >0$. Let $T_i:E\to E, i=1,2,..., r$ be a finite family of $k$-striclty pseudocontractive mappings such that $\Omega:=\underset{i=1}{\overset{r}\bigcap}F(T_i)\ne \emptyset.$ Let $\{a_i\}_{i=1}^r$ be such that
 $0 < a_i<\Bigl(\frac{q k^{q-1}}{d_q}\Bigr)^{\frac{1}{q-1}}, i=1,2,...,r$ and define $T_{a_i}=(1-a_i)I+a_i T_i$ and $S_{a_i}=(1-\omega_i)I+\omega_iT_{a_i}, i=1,2,...,r.$ Let $u\in E$ be fixed and let $\{x_n\}_{n\ge 1}$ be a sequence in $E$ generated iteratively by 
\begin{eqnarray}
x_0\in E, x_{n+1}=\alpha_{n+1}u+(I-\alpha_{n+1}\lambda A)S_{a_{n+1}}x_n, \;n\ge 0,
\end{eqnarray}where $S_{a_n}=S_{{n\;{\rm mod\;}r}}$ and mod function takes values in $\{1,2,...,r\}.$ Then, $\{x_n\}_{n\ge 0}$ converges strongly to a solution of the variational inequality \eqref{Ma}.
\end{cor}

\begin{cor}
Let $E$ be a $q$-uniformly smooth real Banach space which admits weakly sequentially continuous generalized duality mapping, let $A:E\to E$ be an $L$-Lipschitzian strongly accretive mapping with a constant $\eta >0$. Let $T:E\to E,$ be a $k$-striclty pseudocontractive mapping such that $F(T)\ne \emptyset.$ Let be such that
 $0 < a<\Bigl(\frac{q k^{q-1}}{d_q}\Bigr)^{\frac{1}{q-1}}$ and define $T_{a}=(1-a)I+aT.$ Let $u\in E$ be fixed and let $\{x_n\}_{n\ge 1}$ be a sequence in $E$ generated iteratively by 
\begin{eqnarray}
x_0\in E, x_{n+1}=\alpha_{n+1}u+(I-\alpha_{n+1}\lambda A)T_ax_n, \;n\ge 0,
\end{eqnarray} then, $\{x_n\}_{n\ge 0}$ converges strongly to a unique $x^\prime \in F(T)$ which is a solution of the variational inequality \eqref{u}
\end{cor}

\begin{rmk}
\noindent Prototype for our iteration parameter $\{\alpha_n\}_{n\ge 1}$ (see e.g. \cite{27}) is given by 
\begin{eqnarray*}\label{e11}
 \alpha_n=\left\{
\begin{array}{ll}
\frac{1}{\sqrt{n}}\;if \;n\;{\rm is\;odd}, \\
\\
 \frac{1}{\sqrt{n}-1}\;if \;n\;{\rm is\;even}. \\
\end{array}
\right.
\end{eqnarray*}
 \noindent  If we assume that $r$ is odd, since the case $r$ being even is similar, it is not difficult to see that 
 \begin{eqnarray*}
 \frac{\alpha_n}{\alpha_{n+r}}=\left\{
\begin{array}{ll}
 \frac{\sqrt{n+r}-1}{\sqrt{n}}\;if \;n\;{\rm is\;odd}, \\
\\
 \frac{\sqrt{n+r}}{\sqrt{n}-1}\;if \;n\;{\rm is\;even}. \\
\end{array}
\right.
\end{eqnarray*}
\end{rmk}

\begin{rmk}
It is easy to see that Corollary \ref{piu} is an obvious improvement on the corresponding results of Xu \cite{27} in the sense that the condition
$$\Omega= F(T_{r}T_{r-1}...T_{1})=F(T_{1}T_{r-1}...T_{2})=...=F(T_{r-1}...T_{1}T_{r})$$ which was imposed by Xu is dispensed with.
\end{rmk}

\begin{rmk}
When $E=H$, a Hilbert space, $\lambda=1$ and $A$ is a bounded linear strongly positive operator, our iteration process \eqref{Mama} reduces to the iteration scheme studied by Xu \cite{27}. If the fixed vector $u\in E$ is identically equal to the zero vector of $H$, $A$ a strongly monotone operator and we consider a single nonexpansive mapping, \eqref{Mama} reduces to the scheme studied by Yamada \cite{Yamada}. We recall that the results of Xu \cite{27} and  Yamada \cite{Yamada} remain in Hilbert spaces. Our theorems, therefore, extend, generalize, improve and unify the corresponding results of these authors and that of a host of other authors in the more general setting of $q$-uniformly real Banach spaces and for the more general class of strongly accretive mappings. Our corollaries, applications and method of proof are of independent interest.
\end{rmk}

\end{document}